\documentclass[11pt]{article}

\usepackage[utf8]{inputenc}
\usepackage[slovak,english]{babel}

\usepackage{enumerate}

\usepackage{geometry}
\usepackage[leqno]{amsmath}

\usepackage{amssymb}
\usepackage{xfrac}
\usepackage{dsfont}

\usepackage[hypcap=true]{caption}

\usepackage{hyperref}
\usepackage{pdfpages}
\hypersetup{colorlinks=true,linkcolor=blue,citecolor=green!50!black}

\usepackage[thmmarks,thref,hyperref]{ntheorem}
\makeatletter
\newtheoremstyle{nonumberbreakwithoutbrackets}%
  {\item[\rlap{\vbox{\hbox{\hskip\labelsep \theorem@headerfont ##1\theorem@separator}\hbox{\strut}}}]}%
  {\item[\rlap{\vbox{\hbox{\hskip\labelsep \theorem@headerfont ##1\ ##3\theorem@separator}\hbox{\strut}}}]}
  
  \newcommand{\listintertext}{\@ifstar\listintertext@\listintertext@@}
  \newcommand{\listintertext@}[1]{
  	\hspace*{-\@totalleftmargin}#1}
  \newcommand{\listintertext@@}[1]{
  	\hspace{-\leftmargin}#1}
\makeatother

\theoremstyle{nonumberplain}
\theoremseparator{. }
\newtheorem{theorem}{Theorem}

\theoremheaderfont{\normalfont\bfseries}
\theorembodyfont{\upshape}
\theoremstyle{nonumberbreakwithoutbrackets}
\theoremseparator{. }
\theoremsymbol{\ensuremath{\blacksquare}}
\newtheorem{bigproof}{Proof}

\theoremheaderfont{\normalfont}
\theoremsymbol{\ensuremath{\blacksquare}}
\newtheorem{smallproof}{Proof}

\usepackage{pgf,tikz}
\usetikzlibrary{arrows,automata, patterns, calc}

\def\thedate{\today}
\def\me{Igor Fabrici, Jochen Harant, Samuel Mohr, Jens M. Schmidt}
\def\departA{Ilmenau University of Technology, Department of Mathematics, Ilmenau, Germany}
\def\departB{Pavol Jozef Šafárik University, Institute of Mathematics, Košice, Slovakia}
\def\titlevar{Longer Cycles in Essentially 4-Connected Planar Graphs}

\title\titlevar
\author\me
\date\thedate

\hypersetup{pdftitle=\titlevar,pdfsubject=Note,pdfauthor=\me}

\newcommand\newparagraph{\vspace*{4ex}}

\parskip\medskipamount 

\makeatletter
\newcommand\need[1]{\par \penalty-100 \begingroup 
   \dimen@\pagegoal \advance\dimen@-\pagetotal 
   \ifdim #1>\dimen@ 
      \ifdim\dimen@>\z@ \vskip -\pagedepth plus 1fil \fi 
      \break 
   \fi \endgroup}

\let\int\relax
\DeclareMathOperator\int{int}

\let\circ\relax
\DeclareMathOperator\circ{circ}

\let\oldoverline\overline
\renewcommand{\overline}[2][3]{{}\mkern#1mu\oldoverline{\mkern-#1mu#2}}

\newcommand{\ml}{l\kern-0.55mm\char39\kern-0.3mm}

\makeatother

\begin{document}
\begin{center}
{\bf \Large \titlevar}\\[3mm]
{\bf Igor Fabrici\textsuperscript{\textnormal{a,}}\footnote{\label{DAAD}partially supported by DAAD, Germany (as part of BMBF) and by the Ministry of Education, Science, Research and Sport of the Slovak Republic within the project 57320575.}\textsuperscript{\textnormal{,}}\footnote{partially supported by Science and Technology Assistance Agency under the contract
No. APVV-15-0116 and by the Slovak VEGA Grant 1/0368/16.}, Jochen Harant\textsuperscript{\textnormal{b,\ref{DAAD}}}, Samuel Mohr\textsuperscript{\textnormal{b,\ref{DAAD}}}, Jens M.\ Schmidt\textsuperscript{\textnormal{b,\ref{DAAD}}}}\\
\textsuperscript{a} \departB \\
\textsuperscript{b} \departA \\
\end{center}

\hrulefill

\begin{abstract}
A planar 3-connected graph $G$ is called \emph{essentially $4$-connected} if, for every 3-separator $S$, at least one of the two components of $G-S$ is an isolated vertex. Jackson and Wormald proved that the length $\circ(G)$ of a longest cycle of any essentially 4-connected planar graph $G$ on $n$ vertices is at least $\frac{2n+4}{5}$ and Fabrici, Harant and \foreignlanguage{slovak}{Jendroľ} improved this result to $\circ(G)\geq \frac{1}{2}(n+4)$. 
In the present paper, we prove that an essentially 4-connected planar graph on $n$ vertices contains a cycle of length at least $\frac{3}{5}(n+2)$ and that such a cycle can be found in time $O(n^2)$.
\end{abstract}

\bigskip

\noindent {\bf Keywords.} essentially 4-connected planar graph, longest cycle, circumference, shortness coefficient.

\bigskip

For a finite and simple graph $G$ with vertex set $V(G)$ and edge set $E(G)$, let $N(x)$ and $d(x) = |N (x)|$ denote the neighborhood and the degree of any $x\in V(G)$ in $G$, respectively.
The \emph{circumference} $\circ(G)$ of a graph $G$ is the length of a longest cycle of $G$.
A subset $S \subseteq  V(G)$ is an \emph{$s$-separator} of $G$ if $|S| = s$ and $G-S$ is disconnected. From now on, let $G$ be a 3-connected planar graph. 
For every 3-separator $S$ of $G$, it is well-known that $G- S$ has exactly two components. 
We call $S$ \emph{trivial} if at least one component of $G-S$ is a single vertex. If every 3-separator $S$ of $G$ is trivial, we call the 3-connected graph $G$ \emph{essentially $4$-connected}. In the present paper, we are interested in lower bounds on the circumference of essentially 4-connected planar graphs.

Jackson and Wormald~\cite{JW} proved that $\circ(G)\geq \frac{2n+4}{5}$ for every essentially 4-connected planar graph on $n$ vertices and presented an infinite family of essentially 4-connected planar graphs $G$ such that  $\circ(G)\leq c\cdot n$  for each real constant $c>\frac23$. 
Moreover, there is a construction of infinitely many essentially 4-connected planar graphs with $\circ(G)=\frac{2}{3}(n+4)$ (for example see~\cite{FHJ}).
It is open whether there exists an essentially $4$-connected planar graph $G$ on $n$ vertices with $\circ(G)<\frac{2}{3}(n+4)$.
Further results on the length of longest cycles in essentially 4-connected planar graphs can be found in~\cite{FHJ,GM,Z}.

Fabrici, Harant and \foreignlanguage{slovak}{Jendroľ}~\cite{FHJ} extended the result of Jackson and Wormald by proving that $\circ(G)\ge \frac{1}{2}(n+4)$ for every essentially 4-connected planar graph $G$ on $n$ vertices.

Our result is presented in the following theorem. 

\begin{theorem}\label{theorem} 
For any essentially $4$-connected planar graph $G$ on $n$ vertices, $\circ(G)\geq \frac{3}{5}(n+2)$. 
\end{theorem}

We remark that the assertion of \hyperref[theorem]{the theorem} can be improved to $\circ(G)\geq \frac{3}{5}(n+4)$ if $n\geq 16$. This follows from using Lemma~5 in~\cite{FHJ} and a more special version of the forthcoming inequality~(\ref{eq:numberfaces}). We will also show how cycles of $G$ of length at least $\frac{3}{5}(n+2)$ can be found in quadratic time. 
 
\newparagraph

Let $C$ be a plane cycle and let $B$ be a set disjoint from $V(C)$.
A plane graph $H$ is called a \emph{$(B,C)$-graph} if $B\cup V(C)$ is the vertex set of $H$, the cycle $C$ is an induced subgraph of $H$, the subgraph of $H$ induced by $B$ is edgeless, and each vertex of $B$ has degree 3 in $H$.
The vertices in $B$ are called \emph{outer vertices} of $C$. 

A face $f$ of $H$ is called \emph{minor} (\emph{major}) if it is incident with at most one (at least two) outer vertices. 
Note that $f$ is incident with no outer vertex if and only if $C$ is the facial cycle of $f$.

For every $(B,C)$-graph $H$, let $\mu(H)$ denote the number of minor faces of $H$. Then
\begin{align}\label{eq:numberfaces}
\mu(H)\geq |V(H)|-|V(C)|+2. 
\end{align}

\begin{bigproof}[of (\ref{eq:numberfaces})]
Let $H$ be a smallest counterexample. 
Since $B=\emptyset$ implies $|V(H)|=|V(C)|$ and $\mu(H)=2$, which satisfies the inequality~(\ref{eq:numberfaces}), we may assume that $B$ is non-empty. 
For each vertex $y\in B$, the three neighbors of $y$ divide $C$ into three internally disjoint paths $P_1(y)$, $P_2(y)$, and $P_3(y)$ with endvertices in $N(y)$. 
We may assume that $|V(P_1(y))|\leq |V(P_2(y))|\leq |V(P_3(y))|$ and define $\phi(y)=|V(P_1(y))|+|V(P_2(y))|-1$ in this case.

Let $x\in B$ be chosen such that $\phi(x) = \min\{\phi(y)\mid y\in B\}$. 
Consider the two cycles $A_1$ and $A_2$ induced by $V(P_1(x))\cup \{x\}$ and $V(P_2(x))\cup\{x\}$, respectively. 
We claim that the interior of $A_1$ as well as the interior of $A_2$ is a face of $H$ and hence, both are minor faces. 
Suppose that there is a vertex $z$ in the interior of $A_i$ for $i\in\{1,2\}$. 
Then $\phi(z) = |V(P_1(z))| + |V(P_2(z))| - 1 \leq\max\{\,|V(P_1(x))|, |V(P_2(x))|\,\} < |V(P_1(x))| + |V(P_2(x))| - 1 = \phi(x)$, which contradicts the choice of $x$.

Let $H'=H-x$. Note that $H'$ is a $((B\setminus\{x\}),C)$-graph  and has less vertices than $H$. 
Then $|V(H')|=|V(H)|-1$, $\mu(H')\leq \mu(H)-1$, and $\mu(H')\geq |V(H')|-|V(C)|+2$, hence $\mu(H)\ge 1+\mu(H')\ge 1+|V(H')|-|V(C)|+2=|V(H)|-|V(C)|+2$.
\end{bigproof}

\begin{bigproof}[of the Theorem]
Let $G$ be an essentially $4$-connected plane graph on $n$ vertices. 
If $G$ has at most 10 vertices, then it is well known that $G$ is Hamiltonian~\cite{Dillencourt}. In this case, we are done, since $n\geq \frac{3}{5}(n+2)$ for $n\geq 3$. 
Thus, we assume $n\geq 11$. 
A cycle $C$ of $G$ is called an \emph{outer-independent-$3$-cycle} (OI3-cycle) if $V(G)\setminus V(C)$ is an independent set of vertices and $d(x) = 3$ for every $x\in V(G)\setminus V(C)$.
An edge $a = xy\in E(C)$ of a cycle $C$ is called an \emph{extendable edge} of $C$ if $x$ and $y$ have a common neighbor in $V(G)\setminus V(C)$.

In \cite{FHJ}, it is shown that every essentially 4-connected planar graph $G$ on $n\geq 11$ vertices contains an OI3-cycle. In this proof, let $C$ be a longest OI3-cycle of $G$, let $c=|V(C)|$, and let $H$ be the graph obtained from $G$ by removing all \emph{chords} of $C$, i.\,e.\ by removing all edges in $E(G)\setminus E(C)$ that connect vertices of $C$. Clearly, $C$ does not contain an extendable edge.
Obviously, $H$ is a $(B,C)$-graph, with $B=V(H)\setminus V(C)$. 

For the number $\mu$ of minor faces of $H$, we have by (\ref{eq:numberfaces})
\begin{align}
\nonumber
\mu&\geq n-c+2. \\
\intertext{Moreover, we will show}
\label{eq:facestovertices}
6\,\mu&\leq 4\,c
\end{align}
and then, \hyperref[theorem]{the theorem} follows immediately. 

\begin{bigproof}[of (\ref{eq:facestovertices})]

An edge $e$ of $C$ is incident with exactly two faces $f_1$ and $f_2$ of $H$. 
In this case, we say $f_1$ is \emph{opposite} to $f_2$ with respect to $e$. 
A face $f$ of $H$ is called \emph{$j$-face} if it is incident with exactly $j$ edges of $C$ and the edges of $C$ incident with $f$ are called \emph{$C$-edges} of $f$. 
Because $C$ does not contain an extendable edge, we have $j\geq 2$ for every minor $j$-face of $H$.

We define a weight function $w_0$ on the set $F(H)$ of faces of $H$, by setting weight $w_0(f)=6$ for every minor face $f$ of $H$ and weight $w_0(f)=0$ for every major face $f$ of $H$. Then $\sum_{f\in F(H)}w_0(f)= 6\,\mu$.
Next, we redistribute the weights of faces of $H$ by the rules \textbf{R1} and  \textbf{R2}.

\begin{description}
\item[Rule R1.] A minor $2$-face $f$ of $H$ sends weight $1$ through both $C$-edges to the opposite (possibly identical) faces.
\item[Rule R2.] A minor $3$-face $f$ of $H$ with $C$-edges $ux$, $xy$, and $yz$ sends weight $1$ through its \emph{middle $C$-edge} $xy$ to the opposite face.
\end{description}
Let $w_1$ denote the new weight function; clearly, $\sum_{f\in F(H)}w_1(f)=6\,\mu$ still holds. 

For the proof of (\ref{eq:facestovertices}), we will show
\begin{align}
w_1(f)\leq 2\,j \text{ for each $j$-face $f$ of $H$.}\label{eq:aftermoving}
\end{align}
To see that (\ref{eq:facestovertices}) is a consequence of (\ref{eq:aftermoving}), let each $j$-face $f$ of $H$ satisfying $j\geq 1$ send the weight $\frac{w_1(f)}{j}$ to each of its $C$-edges. 
Note that each $0$-face $f$ is major, thus $w_1(f)=0$.
Hence, the total weight of all minor and major faces is moved to
the edges of $C$. 
Since every edge of $C$ gets weight at most $4$, we obtain $6\,\mu=\sum_{f\in F(H)}w_1(f)\leq 4\,c$, and (\ref{eq:facestovertices}) follows.

\newparagraph

\begin{smallproof}[of (\ref{eq:aftermoving})]
Next we distinguish several cases. In most of them, we construct a cycle $\tilde{C}$ that is obtained from $C$ by replacing a subpath of $C$ with another path. In every case, $\tilde{C}$ will be an OI3-cycle of $G$ that is longer than $C$. This contradicts the choice of $C$ and therefore shows that the considered case cannot occur. Note that all vertices of $C$ in the following figures are different, because  the length of the longest OI3-cycle $C$ in a planar graph on $n\geq 11$ vertices is at least $8$~\cite[Lemma~4(ii)]{FHJ}.

\begin{description}
\need{1cm}
\item[Case 1.] 
\textit{$f$ is a major $j$-face. }

Because $w_0(f)=0$ and $f$ gets weight $\leq 1$ through each of its $C$-edge, we have $w_1(f)\leq j$.

\need{1cm}
\item[Case 2.] 
\textit{$f$ is a minor $2$-face (see Figure~$\ref{fig:2face}$). }

We will show that $f$ does not get any new weight by \textbf{R1} or by \textbf{R2}; this implies $w_1(f)=w_0(f)-(1+1)=4$. Let $xy$ and $yz$ be the $C$-edges of $f$ and $a$ be the outer vertex incident with $f$ (see Figure~$\ref{fig:2face}$).

\begin{figure}[ht]
\begin{center}
\begin{tikzpicture}[scale=0.25, -, 
vertex/.style={circle,fill=black,draw,minimum size=5pt,inner sep=0pt}]

\node[vertex,label=270:{ $x$}] (x) at (0,0) [] {$ $};

\node[vertex,label=270:{ $y$}] (y) at (4,0) [] {$ $};

\node[vertex,label=270:{ $z$}] (z) at (8,0) [] {$ $};

\node[vertex,label=90:{ $a$}] (a) at (4,7) [] {$ $};

\draw (4,4) node {$f$};

\node (p1) at(0,0) [] {$ $};

\node (p2) at(8,0) [] {$ $};

\node (q1) at(-4,1) [] {$ $};

\node[label=0:{ $C$}] (q2) at(12,1) [] {$ $};

\path[-,out=-30, in=180] (q1) edge (p1);

\draw (x) -- (z);

\path[-,out=0, in=-150] (p2) edge (q2);

\draw (x) -- (a);

\draw (z) -- (a);
\end{tikzpicture}
\end{center}
\caption{}\label{fig:2face}
\end{figure}

If $f$ gets new weight by \textbf{R1} or by \textbf{R2} from a face $f'$ opposite to $f$ with respect to a $C$-edge of $f$, then $f'$ is a minor $2$-face or a minor $3$-face of $H$. Without loss of generality, we may assume  that  $f'$ is opposite to $f$ with respect to the edge $yz$. Then $yz$ is a common $C$-edge of $f$ and $f'$ and we distinguish the following subcases.

\begin{description}
\need{1cm}
\item[Case 2a.]
\textit{$f'$ is a $2$-face and $xy$ is a $C$-edge of $f'$. }

Then $\{x,z\}$ is the neighborhood of $y$ in $G$, which contradicts the $3$-connectedness of $G$.

\need{1cm}
\item[Case 2b.]
\textit{$f'$ is a $2$-face and $xy$ is not a $C$-edge of $f'$ (see Figure~$\ref{fig:2face2}$). }

Then a longer OI3-cycle $\tilde{C}$ is obtained from $C$ by replacing the path $(x,y,z,u)$ with the path $(x,a,z,y,b,u)$, which gives a contradiction.

\begin{figure}[ht]
\begin{minipage}[b]{.5\textwidth}
\begin{center}
\begin{tikzpicture}[scale=0.25, -, 
vertex/.style={circle,fill=black,draw,minimum size=5pt,inner sep=0pt}]

\node[vertex,label=270:{ $x$}] (x) at (0,0) [] {$ $};

\node[vertex,label=270:{ $y$}] (y) at (4,0) [] {$ $};

\node[vertex,label=270:{ $z$}] (z) at (8,0) [] {$ $};

\node[vertex,label=270:{ $u$}] (u) at (12,0) [] {$ $};

\node[vertex,label=90:{ $a$}] (a) at (4,7) [] {$ $};

\draw (4,4) node {$f$};

\node[vertex,label=-90:{ $b$}] (b) at (8,-7) [] {$ $};

\node (p1) at(0,0) [] {$ $};

\node (p2) at(12,0) [] {$ $};

\node (q1) at(-4,1) [] {$ $};

\node[label=0:{ $C$}] (q2) at (16,1) [] {$ $};

\path[-,out=-30, in=180] (q1) edge (p1);

\draw (x) -- (u);

\path[-,out=0, in=-150] (p2) edge (q2);

\draw (x) -- (a);

\draw (z) -- (a);

\draw (y) -- (b);

\draw (u) -- (b);
\end{tikzpicture}
\end{center}
\caption{}\label{fig:2face2}
\end{minipage}\hfil %
\begin{minipage}[b]{.5\textwidth}
\begin{center}

\begin{tikzpicture}[scale=0.25, -, 
vertex/.style={circle,fill=black,draw,minimum size=5pt,inner sep=0pt}]

\node[vertex,label=270:{ $x$}] (x) at (0,0) [] {$ $};

\node[vertex,label=270:{ $y$}] (y) at (4,0) [] {$ $};

\node[vertex,label=270:{ $z$}] (z) at (8,0) [] {$ $};

\node[vertex,label=270:{ $u$}] (u) at (12,0) [] {$ $};

\node[vertex,label=90:{ $a$}] (a) at (4,7) [] {$ $};

\draw (4,4) node {$f$};

\node[vertex,label=-90:{ $b$}] (b) at (6,-7) [] {$ $};

\node (p1) at(0,0) [] {$ $};

\node (p2) at(12,0) [] {$ $};

\node (q1) at(-4,1) [] {$ $};

\node[label=0:{ $C$}] (q2) at(16,1) [] {$ $};

\path[-,out=-30, in=180] (q1) edge (p1);

\draw (x) -- (u);

\path[-,out=0, in=-150] (p2) edge (q2);

\path[dashed,out=-45, in=-135] (y) edge (u);

\draw (x) -- (a);

\draw (z) -- (a);

\draw (x) -- (b);

\draw (u) -- (b);

\end{tikzpicture}
\end{center}
\caption{}\label{fig:2face3}
\end{minipage}
\end{figure}

\need{1cm}
\item[Case 2c.]
\textit{$f'$ is a $3$-face. }

Since $f'$ sends weight to $f$, then, by rule \textbf{R2}, a $C$-edge of $f$ is the middle $C$-edge of $f'$. It follows that both $C$-edges of $f$ are also $C$-edges of $f'$ and the situation as shown in Figure~\ref{fig:2face3} occurs. The edge $yu$ exists in $G$, because otherwise $d(y)=2$ and $G$ would not be 3-connected. Then $\tilde{C}$ is obtained  from $C$ by replacing the path $(x,y,z,u)$ with the path $(x,a,z,y,u)$.
\end{description}

\need{1cm}
\item[Case 3.] 
\textit{$f$ is a minor $3$-face (see Figure~$\ref{fig:3face}$). }

Since $f$ looses weight $1$ by rule \textbf{R2} and possible gets weight $w$ by \textbf{R1} or by \textbf{R2}, we have $w_1(f')=5+w$.

If $w\leq 1$, then we are done.

\begin{figure}[ht]
\begin{center}
\begin{tikzpicture}[scale=0.25, -, 
vertex/.style={circle,fill=black,draw,minimum size=5pt,inner sep=0pt}]

\node[vertex,label=270:{ $v$}] (v) at (-4,0) [] {$ $};

\node[vertex,label=270:{ $x$}] (x) at (0,0) [] {$ $};

\node[vertex,label=270:{ $y$}] (y) at (4,0) [] {$ $};

\node[vertex,label=270:{ $z$}] (z) at (8,0) [] {$ $};

\node[vertex,label=90:{ $a$}] (a) at (2,7) [] {$ $};

\draw (2,4) node {$f$};

\node (p1) at(-4,0) [] {$ $};

\node (p2) at(8,0) [] {$ $};

\node (q1) at(-8,1) [] {$ $};

\node[label=0:{ $C$}] (q2) at(12,1) [] {$ $};

\path[-,out=-30, in=180] (q1) edge (p1);

\draw (v) -- (z);

\path[-,out=0, in=-150] (p2) edge (q2);

\draw (v) -- (a);

\draw (z) -- (a);
\end{tikzpicture}
\end{center}
\caption{}\label{fig:3face}
\end{figure}

If $w\geq 2$, then $f$ do not get any weight through the edge $xy$ from the opposite face $f'$. 
Otherwise, if $f'$ is a $2$-face, then we have the situation as in Case 2c and if $f'$ is a $3$-face, then $w=1$, with contradiction in both cases. 
Hence, $f$ gets weight 1 through $vx$ from the opposite face $f_1$ and weight 1 through $yz$ from the opposite face $f_2$. 
Clearly, $f_1\neq f_2$ and they are not simultaneously $3$-faces.

\begin{description}
\need{1cm}
\item[Case 3a.]
\textit{Both $f_1$ and $f_2$ are $2$-faces. }

Then  the situation is as illustrated in Figure~\ref{fig:3face22} and $\tilde{C}$ is obtained from $C$ by replacing the path $(w,v,x,y,z,u)$ with the path $(w,b,x,v,a,z,y,c,u)$. 
Note that $b\neq c$, because $d(b)=d(c)=3$. 

\begin{figure}[ht]
\begin{minipage}[b]{.5\textwidth}
\begin{center}
\begin{tikzpicture}[scale=0.25, -, 
vertex/.style={circle,fill=black,draw,minimum size=5pt,inner sep=0pt}]

\node[vertex,label=270:{ $w$}] (w) at (-8,0) [] {$ $};

\node[vertex,label=270:{ $v$}] (v) at (-4,0) [] {$ $};

\node[vertex,label=270:{ $x$}] (x) at (0,0) [] {$ $};

\node[vertex,label=270:{ $y$}] (y) at (4,0) [] {$ $};

\node[vertex,label=270:{ $z$}] (z) at (8,0) [] {$ $};

\node[vertex,label=270:{ $u$}] (u) at (12,0) [] {$ $};

\node[vertex,label=90:{ $a$}] (a) at (2,7) [] {$ $};

\draw (2,4) node {$f$};

\node[vertex,label=-90:{ $b$}] (b) at (-4,-7) [] {$ $};

\node[vertex,label=-90:{ $c$}] (c) at (8,-7) [] {$ $};

\node (p1) at(-8,0) [] {$ $};

\node (p2) at(12,0) [] {$ $};

\node (q1) at(-12,1) [] {$ $};

\node[label=0:{ $C$}] (q2) at(16,1) [] {$ $};

\path[-,out=-30, in=180] (q1) edge (p1);

\draw (w) -- (u);

\path[-,out=0, in=-150] (p2) edge (q2);

\draw (v) -- (a);

\draw (z) -- (a);

\draw (w) -- (b);

\draw (x) -- (b);

\draw (y) -- (c);

\draw (u) -- (c);

\end{tikzpicture}
\caption{}\label{fig:3face22}
\end{center}
\end{minipage}\hfil %
\begin{minipage}[b]{.5\textwidth}
\begin{center}

\begin{tikzpicture}[scale=0.25, -, 
vertex/.style={circle,fill=black,draw,minimum size=5pt,inner sep=0pt}]

\node[vertex,label=270:{ $w$}] (w) at (-8,0) [] {$ $};

\node[vertex,label=270:{ $v$}] (v) at (-4,0) [] {$ $};

\node[vertex,label=270:{ $x$}] (x) at (0,0) [] {$ $};

\node[vertex,label=270:{ $y$}] (y) at (4,0) [] {$ $};

\node[vertex,label=270:{ $z$}] (z) at (8,0) [] {$ $};

\node[vertex,label=270:{ $u$}] (u) at (12,0) [] {$ $};

\node[vertex,label=90:{ $a$}] (a) at (2,7) [] {$ $};

\draw (2,4) node {$f$};

\node[vertex,label=-90:{ $b$}] (b) at (-4,-7) [] {$ $};

\node[vertex,label=-90:{ $c$}] (c) at (6,-7) [] {$ $};

\node (p1) at(-8,0) [] {$ $};

\node (p2) at(12,0) [] {$ $};

\node (q1) at(-12,1) [] {$ $};

\node[label=0:{ $C$}] (q2) at(16,1) [] {$ $};

\path[-,out=-30, in=180] (q1) edge (p1);

\draw (w) -- (u);

\path[-,out=0, in=-150] (p2) edge (q2);

\draw (v) -- (a);

\draw (z) -- (a);

\draw (w) -- (b);

\draw (x) -- (b);

\draw (x) -- (c);

\draw (u) -- (c);

\end{tikzpicture}
\caption{}\label{fig:3face23}
\end{center}
\end{minipage}
\end{figure}

\need{1cm}
\item[Case 3b.]
\textit{$f_1$ is a $2$-face and $f_2$ is a $3$-face. }

Then $e_2=yz$ is the middle $C$-edge of $f_2$, as shown in Figure~\ref{fig:3face23}, and $\tilde{C}$ is obtained from $C$ by replacing the path $(w,v,x,y,z,u)$ with the path $(w,v,a,z,y,x,c,u)$.
\end{description}

\need{1cm}
\item[Case 4.] 
\textit{$f$ is a minor $4$-face (see Figure~$\ref{fig:4face}$). }

\begin{figure}[ht]
\begin{center}
\begin{tikzpicture}[scale=0.25, -, 
vertex/.style={circle,fill=black,draw,minimum size=5pt,inner sep=0pt}]

\node[vertex,label=270:{ $v$}] (v) at (-8,0) [] {$ $};

\node[vertex,label=270:{ $w$}] (w) at (-4,0) [] {$ $};

\node[vertex,label=270:{ $x$}] (x) at (0,0) [] {$ $};

\node[vertex,label=270:{ $y$}] (y) at (4,0) [] {$ $};

\node[vertex,label=270:{ $z$}] (z) at (8,0) [] {$ $};

\node[vertex,label=90:{ $a$}] (a) at (0,7) [] {$ $};

\draw (0,4) node {$f$};

\node (p1) at(-8,0) [] {$ $};

\node (p2) at(8,0) [] {$ $};

\node (q1) at(-12,1) [] {$ $};

\node[label=0:{ $C$}] (q2) at(12,1) [] {$ $};

\path[-,out=-30, in=180] (q1) edge (p1);

\draw (v) -- (z);

\path[-,out=0, in=-150] (p2) edge (q2);

\draw (v) -- (a);

\draw (z) -- (a);
\end{tikzpicture}
\end{center}
\caption{}\label{fig:4face}
\end{figure}

If $w_1(f)=w_0(f)+w=6+w$ and $w\leq 2$, then we are done. 

If otherwise $w\geq 3$, there are at least three edges $e_1$, $e_2$, and $e_3$ among the four $C$-edges $vw$, $wx$, $xy$, and $yz$ of $f$ such that $f$ gets weight from  minor faces which are opposite to $f$ with respect to $e_1$, $e_2$, and $e_3$, respectively. 

\begin{description}
\need{1cm}
\item[Case 4a.]
\textit{$w=3$ and $\{e_1,e_2,e_3\}=\{vw,wx,xy\}$. }

Then no edge of $\{e_1,e_2,e_3\}$ is the middle $C$-edge of a minor $3$-face and $yz$ is not a $C$-edge of a minor $2$-face. 
We have the situation of Figure~\ref{fig:4face22} and one of the edges $vx$ or $xz$ exists in  $G$, because otherwise $x$ would have degree 2 in $G$.

Then $\tilde{C}$ is obtained again from $C$ by replacing the path $(v,w,x,y,z)$ with the path $(v,x,w,c,y,z)$ or with the path $(v,w,c,y,x,z)$, respectively.

\begin{figure}[ht]
\begin{minipage}[b]{.5\textwidth}
\begin{center}
\begin{tikzpicture}[scale=0.25, -, 
vertex/.style={circle,fill=black,draw,minimum size=5pt,inner sep=0pt}]

\node[vertex,label=270:{ }] (f) at (-12,0) [] {$ $};

\node[vertex,label=270:{ $v$}] (v) at (-8,0) [] {$ $};

\node[vertex,label=270:{ $w$}] (w) at (-4,0) [] {$ $};

\node[vertex,label=270:{ $x$}] (x) at (0,0) [] {$ $};

\node[vertex,label=270:{ $y$}] (y) at (4,0) [] {$ $};

\node[vertex,label=270:{ $z$}] (z) at (8,0) [] {$ $};

\node[vertex,label=90:{ $a$}] (a) at (0,7) [] {$ $};

\draw (0,4) node {$f$};

\node[vertex,label=-90:{ $b$}] (b) at (-8,-7) [] {$ $};

\node[vertex,label=-90:{ $c$}] (c) at (0,-7) [] {$ $};

\node (p1) at(-12,0) [] {$ $};

\node (p2) at(8,0) [] {$ $};

\node (q1) at(-16,1) [] {$ $};

\node[label=0:{ $C$}] (q2) at(12,1) [] {$ $};

\path[-,out=-30, in=180] (q1) edge (p1);

\draw (f) -- (z);

\path[-,out=0, in=-150] (p2) edge (q2);

\path[dashed,out=45, in=135] (v) edge (x);

\path[dashed,out=45, in=135] (x) edge (z);

\draw (v) -- (a);

\draw (z) -- (a);

\draw (f) -- (b);

\draw (w) -- (b);

\draw (w) -- (c);

\draw (y) -- (c);
\end{tikzpicture}
\caption{}\label{fig:4face22}
\end{center}
\end{minipage}\hfil %
\begin{minipage}[b]{.5\textwidth}
\begin{center}
\begin{tikzpicture}[scale=0.25, -, 
vertex/.style={circle,fill=black,draw,minimum size=5pt,inner sep=0pt}]

\node[vertex,label=270:{ $t$}] (f) at (-12,0) [] {$ $};

\node[vertex,label=270:{ $v$}] (v) at (-8,0) [] {$ $};

\node[vertex,label=270:{ $w$}] (w) at (-4,0) [] {$ $};

\node[vertex,label=270:{ $x$}] (x) at (0,0) [] {$ $};

\node[vertex,label=270:{ $y$}] (y) at (4,0) [] {$ $};

\node[vertex,label=270:{ $z$}] (z) at (8,0) [] {$ $};

\node[vertex,label=90:{ $a$}] (a) at (0,7) [] {$ $};

\draw (0,4) node {$f$};

\node[vertex,label=-90:{ $b$}] (b) at (-8,-7) [] {$ $};

\node[vertex,label=-90:{ $c$}] (c) at (4,-7) [] {$ $};

\node (p1) at(-12,0) [] {$ $};

\node (p2) at(8,0) [] {$ $};

\node (q1) at(-16,1) [] {$ $};

\node[label=0:{ $C$}] (q2) at(12,1) [] {$ $};

\path[-,out=-30, in=180] (q1) edge (p1);

\draw (f) -- (z);

\path[-,out=0, in=-150] (p2) edge (q2);

\path[dashed,out=45, in=135] (v) edge (y);

\path[dashed,out=45, in=135] (w) edge (y);

\draw (v) -- (a);

\draw (z) -- (a);

\draw (f) -- (b);

\draw (w) -- (b);

\draw (x) -- (c);

\draw (z) -- (c);
\end{tikzpicture}
\caption{}\label{fig:4face2f2}
\end{center}
\end{minipage}
\end{figure}

\need{1cm}
\item[Case 4b.]
\textit{$w=3$, $\{e_1,e_2,e_3\}=\{vw,xy,yz\}$ and $wx$ is not a $C$-edge of a minor $3$-face. }

Then  $vw$ is not the middle $C$-edge of a minor  $3$-face opposite to $f$. We have the situation of Figure~\ref{fig:4face2f2} and one of the edges $vy$ or $wy$ exists in  $G$, because otherwise $y$ would have degree 2 in $G$.

Note that $b\neq c$, because $d(b)=d(c)=3$. Then $\tilde{C}$ is obtained from $C$ by replacing the path $(t,v,w,x,y,z)$ with the path $(t,b,w,v,y,x,c,z)$ or with the path $(t,v,w,y,x,c,z)$.

\need{1cm}
\item[Case 4c.]
\textit{$w=3$, $\{e_1,e_2,e_3\}=\{vw,xy,yz\}$ and $wx$ is  a $C$-edge of a minor $3$-face. }

Then $vw$ is the middle $C$-edge of a minor $3$-face opposite to $f$ (see Figure~\ref{fig:4face32}).

Then at least one of the edges $vy$ or $wy$ exists, because otherwise $y$ would have degree 2 in $G$, and $\tilde{C}$ is obtained from $C$ by replacing the path $(t,v,w,x,y,z)$ with the path $(t,b,x,w,v,y,z)$ or with the path $(t,v,w,y,x,c,z)$.

\begin{figure}[ht]
\begin{minipage}[b]{.5\textwidth}
\begin{center}
\begin{tikzpicture}[scale=0.25, -, 
vertex/.style={circle,fill=black,draw,minimum size=5pt,inner sep=0pt}]

\node[vertex,label=270:{ $t$}] (v) at (-8,0) [] {$ $};

\node[vertex,label=270:{ $v$}] (w) at (-4,0) [] {$ $};

\node[vertex,label=270:{ $w$}] (x) at (0,0) [] {$ $};

\node[vertex,label=270:{ $x$}] (y) at (4,0) [] {$ $};

\node[vertex,label=270:{ $y$}] (z) at (8,0) [] {$ $};

\node[vertex,label=270:{ $z$}] (t) at (12,0) [] {$ $};

\node[vertex,label=90:{ $a$}] (a) at (4,7) [] {$ $};

\draw (4,4) node {$f$};

\node[vertex,label=-90:{ $b$}] (b) at (-2,-7) [] {$ $};

\node[vertex,label=-90:{ $c$}] (c) at (8,-7) [] {$ $};

\node (p1) at(-8,0) [] {$ $};

\node (p2) at(12,0) [] {$ $};

\node (q1) at(-12,1) [] {$ $};

\node[label=0:{ $C$}] (q2) at(16,1) [] {$ $};

\path[-,out=-30, in=180] (q1) edge (p1);

\draw (v) -- (t);

\path[-,out=0, in=-150] (p2) edge (q2);

\path[dashed,out=45, in=135] (w) edge (z);

\path[dashed,out=45, in=135] (x) edge (z);

\draw (w) -- (a);

\draw (t) -- (a);

\draw (v) -- (b);

\draw (y) -- (b);

\draw (y) -- (c);

\draw (t) -- (c);
\end{tikzpicture}
\caption{}\label{fig:4face32}
\end{center}
\end{minipage}\hfil %
\begin{minipage}[b]{.5\textwidth}
\begin{center}
\begin{tikzpicture}[scale=0.25, -, 
vertex/.style={circle,fill=black,draw,minimum size=5pt,inner sep=0pt}]

\node[vertex,label=270:{ $v$}] (v) at (-8,0) [] {$ $};

\node[vertex,label=270:{ $w$}] (w) at (-4,0) [] {$ $};

\node[vertex,label=270:{ $x$}] (x) at (0,0) [] {$ $};

\node[vertex,label=270:{ $y$}] (y) at (4,0) [] {$ $};

\node[vertex,label=270:{ $z$}] (z) at (8,0) [] {$ $};

\node[vertex,label=90:{ $a$}] (a) at (0,7) [] {$ $};

\draw (0,4) node {$f$};

\node[vertex,label=-90:{ $b$}] (b) at (-4,-7) [] {$ $};

\node[vertex,label=-90:{ $c$}] (c) at (4,-7) [] {$ $};

\node (p1) at(-8,0) [] {$ $};

\node (p2) at(8,0) [] {$ $};

\node (q1) at(-12,1) [] {$ $};

\node[label=0:{ $C$}] (q2) at(12,1) [] {$ $};

\path[-,out=-30, in=180] (q1) edge (p1);

\draw (v) -- (z);

\path[-,out=0, in=-150] (p2) edge (q2);

\path[dashed,out=45, in=135] (w) edge (y);

\draw (v) -- (a);

\draw (z) -- (a);

\draw (v) -- (b);

\draw (x) -- (b);

\draw (x) -- (c);

\draw (z) -- (c);
\end{tikzpicture}
\caption{}\label{fig:4face422}
\end{center}
\end{minipage}
\end{figure}

\need{1cm}
\item[Case 4d.]
\textit{$w=4$. }

Then the edges  $vw$, $wx$, $xy$, and $yz$ are $C$-edges of minor $2$-faces of $H$. 
Either a situation similar to Case~4a occurs, a contradiction, or the situation of Figure~\ref{fig:4face422} follows.

Then the edge $wy$ exists in  $G$, because otherwise $d(w)=2$ or $d(y)=2$ in $G$, and $\tilde{C}$ is obtained from $C$ by replacing the path $(v,w,x,y,z)$ with the path $(v,w,y,x,c,z)$.
\end{description}

\need{1cm}
\item[Case 5.] 
\textit{$f$ is a minor $5$-face. }

Let $w_1(f)=w_0(f)+w=6+w$.
If $w\leq 4$, then $w_1(f)\leq 10$ and we are done. 
If $w=5$, then all five $C$-edges of $f$ are also $C$-edges of minor $2$-faces and we have the situation of Figure~\ref{fig:5face222}.

\begin{figure}[ht]
\begin{center}
\begin{tikzpicture}[scale=0.25, -, 
vertex/.style={circle,fill=black,draw,minimum size=5pt,inner sep=0pt}]

\node[vertex,label=270:{ $s$}] (s) at (-12,0) [] {$ $};

\node[vertex,label=270:{ $v$}] (v) at (-8,0) [] {$ $};

\node[vertex,label=270:{ $w$}] (w) at (-4,0) [] {$ $};

\node[vertex,label=270:{ $x$}] (x) at (0,0) [] {$ $};

\node[vertex,label=270:{ $y$}] (y) at (4,0) [] {$ $};

\node[vertex,label=270:{ $z$}] (z) at (8,0) [] {$ $};

\node[vertex,label=270:{ $t$}] (t) at (12,0) [] {$ $};

\node[vertex,label=90:{ $a$}] (a) at (-2,9) [] {$ $};

\draw (-2,5) node {$f$};

\node[vertex,label=-90:{ $b$}] (b) at (-8,-7) [] {$ $};

\node[vertex,label=-90:{ $c$}] (c) at (0,-7) [] {$ $};

\node[vertex,label=-90:{ $d$}] (d) at (8,-7) [] {$ $};

\node (p1) at(-12,0) [] {$ $};

\node (p2) at(12,0) [] {$ $};

\node (q1) at(-16,1) [] {$ $};

\node[label=0:{ $C$}] (q2) at(16,1) [] {$ $};

\path[-,out=-30, in=180] (q1) edge (p1);

\draw (s) -- (t);

\path[-,out=0, in=-150] (p2) edge (q2);

\path[densely dotted,out=45, in=135] (v) edge (x);

\path[dashed,out=45, in=135] (v) edge (z);
\path[dashed,out=45, in=135] (x) edge (z);

\draw (s) -- (a);

\draw (z) -- (a);

\draw (s) -- (b);

\draw (w) -- (b);

\draw (w) -- (c);

\draw (y) -- (c);

\draw (y) -- (d);

\draw (t) -- (d);
\end{tikzpicture}
\end{center}
\caption{}\label{fig:5face222}
\end{figure}

If the edge $vx$  exists, then $\tilde{C}$ is obtained from $C$ by replacing the path $(s,v,w,x)$ with the path $(s,b,w,v,x)$. 

If $vx$ does not exist, then, because $d(v)\geq 3$, $y$ or $z$ is a neighbor of $v$. If the edge $vy$ exists, we get $d(x)=2$, a contradiction. Hence, $vz$ exists and, since $d(x)\geq 3$, $xz$ exists as well. In this case, $\tilde{C}$ is obtained from $C$ by replacing the path $(w,x,y,z)$ with the path $(w,c,y,x,z)$.

\listintertext{The remaining case completes the proof of (\ref{eq:aftermoving}) and therefore the proof of (\ref{eq:facestovertices}). }

\need{1cm}
\item[Case 6.] 
\textit{$f$ is a minor $j$-face with $j\geq 6$. }

Then $w_1(f)=w_0(f)+w=6+w\leq 6+j\leq 2\,j$.
\end{description}
\end{smallproof}
\end{bigproof}
\end{bigproof}

\vspace{-4.8ex}
\paragraph{Algorithm.}
We now show that a cycle of length at least $\frac{3}{5}(n+2)$
in any essentially 4-connected planar graph $G$ on $n$ vertices can be computed in time $O(n^2)$. For $n \leq 10$, we may compute even a longest cycle in constant time, so assume $n \geq 11$.
The existential proof of \hyperref[theorem]{the theorem} proceeds by using a longest not extendable OI3-cycle of $G$. However, it is straight-forward to observe that the proof is still valid when we replace this cycle by an OI3-cycle $C$ that is not extendable and for which none of the local replacements described in the Cases~1--6 can be applied to increase its length (as argued, all these replacements preserve an OI3-cycle).

It suffices to describe how such a cycle $C$ can be computed efficiently; the desired length of $C$ is then implied by \hyperref[theorem]{the theorem}. In~\cite[Lemma~3]{FHJ}, an OI3-cycle of $G$ is obtained by constructing a special Tutte cycle with the aid of Sander's result on Tutte paths~\cite{Sanders}. Using the recent result in~\cite{TuttePaths}, we can compute such Tutte paths and, by prescribing its end vertices accordingly, also the desired Tutte cycle in time $O(n^2)$. This gives an OI3-cycle $C_i$ of $G$.

If $C_i$ is extendable, we compute an extendable edge of $C_i$ and extend $C_i$ to a longer cycle $C_{i+1}$; this takes time $O(n)$ and preserves that $C_{i+1}$ is an OI3-cycle. Otherwise, if there is no extendable edge of $C_i$ (in this case, the length of $C_i$ is at least $8$ due to $n\geq 11$ and~\cite[Lemma~4(ii)]{FHJ}),
we decide in time $O(n)$ whether one of the local replacements of the Cases~1--6 can be applied to $C_i$. 
If so, we apply any such case and obtain the longer OI3-cycle $C_{i+1}$ (which however may be extendable); since all replacements modify only subgraphs of constant size, this can be done in constant time. Iterating this implies a total running time of $O(n^2)$, as the length of the cycle is increased at most $O(n)$ times.

\newparagraph

\end{document}